\documentclass[11pt]{amsart}
\usepackage{graphicx}
\usepackage{amssymb}
\begin{document}

\title{Functionals on the Spaces of $n$-Dimensional Convex Bodies}

\author{Chuanming Zong}

\address{School of Mathematical Sciences, Peking University, Beijing 100871, P. R. China.}
\email{cmzong@math.pku.edu.cn}

\thanks{{\em 2010 Mathematics Subject Classification}: primary 52C17; secondary 52B10, 52C07}

\maketitle

\centerline{\it Dedicated to Professor Wentsun Wu on the occasion of his 95th birthday}

\medskip
\begin{abstract}
In geometry, there are several challenging problems studying numbers associated to convex bodies.
For example, the packing density problem, the kissing number problem, the covering density problem, the packing-covering constant problem, Hadwiger's covering conjecture and Borsuk's partition conjecture.
They are fundamental and fascinating problems about the same objects. However, up to now, both the methodology and the technique applied to them are essentially different. Therefore, a common foundation for them has been much expected. By treating problems of these types as functionals defined on the spaces of $n$-dimensional convex bodies, this paper
tries to create such a foundation. This article suggests an ideal theoretic structure and a couple of research topics such as supderivatives and integral sums of these functionals which seem to be important. In addition, it proves an inequality between the Hausdorff metric and the Banach-Mazur metric and obtains some estimations on the supderivatives.
\end{abstract}

\medskip
\section*{1. Introduction}

\medskip
Let $K$ denote an $n$-dimensional convex body, a convex
and compact set with nonempty interior in $\mathbb{E}^n$, and let $C$ denote
a centrally symmetric one. In particular, let $B^n$ denote the
$n$-dimensional unit ball centered at the origin of $\mathbb{E}^n$ and let
$I^n$ denote the $n$-dimensional unit cube defined by $\{ (x_1, x_2,
\ldots , x_n):\ |x_i|\le {1\over 2}\}$. There are
several important numbers defined on a convex body $K$ such as the
maximal congruent packing density $\delta^c(K)$, the maximal translative packing density $\delta^t(K)$, the maximal
lattice packing density $\delta^l(K)$, the minimal congruent covering density $\theta^c(K)$, the minimal translative
covering density $\theta^t(K)$, the minimal lattice covering density
$\theta^l (K)$, its congruent kissing number $\tau^c(K)$, translative kissing number $\tau^t (K)$, lattice kissing
number $\tau^l(K)$, Hadwiger's covering number $h(K)$ and Borsuk's partition number $b(K)$. For the definitions and history of these numbers, we refer to \cite{bras05}, \cite{conw03}, \cite{feje93} and \cite{zong99}.

\medskip
In 1611, Kepler studied the densities of ball packings and
made the following conjecture:

\smallskip\noindent {\bf Kepler's conjecture.} {\it The maximal
packing density of three-dimensional unit balls is $\pi
/\sqrt{18}$.}

\smallskip
This conjecture and its generalizations have been studied
by many prominent mathematicians, including Gauss, Lagrange, Hermite, Hilbert, Minkowski and others. In 1840, Gauss \cite{gaus40} studied the lattice case and proved
$$\delta^l(B^3)={\pi \over \sqrt{18}}.$$
In 1900, in the third part of his 18th problem, Hilbert
\cite{hilb00} generalized Kepler's conjecture to the following
problem:

\smallskip\noindent {\bf Hilbert's problem.} {\it To determine the maximal packing density of
a given geometric object, for example the unit ball or the regular
tetrahedron.}

\smallskip
In 1904, Minkowski \cite{mink04} discovered a criterion for the
densest lattice packings of a three-dimensional convex body $K$ and
applied it to tetrahedron and octahedron, respectively.
Unfortunately, he made a mistake in the tetrahedron case. In 2000,
based on Minkowski's criterion, Betke and Henk \cite{betk00}
developed an Algorithm by which one can determine the value of
$\delta^l(P)$ for any given three-dimensional polytope $P$.

From 1997 to 2005, in a series of complicated papers (with the
assistance of a computer) Hales published a proof for Kepler's
conjecture (see \cite{hale05} and its references). That is
$$\delta^t(B^3)={\pi \over \sqrt{18}}.$$

Let $\mathcal{K}^n$ denote the family of all $n$-dimensional convex
bodies and let $\mathcal{C}^n$ denote the family of all
$n$-dimensional centrally symmetric convex bodies. By the
definitions of $\delta^c(K)$, $\delta^t(K)$ and $\delta^l(K)$ it is easy to see that
$$\delta^l (K)\le \delta^t(K)\le \delta^c(K)\le 1$$
holds for all $K\in \mathcal{K}^n$, and where the equalities hold when
$K$ is an $n$-dimensional parallelopiped. On the other hand, it is natural to seek the
optimal lower bounds for these numbers. We define
$$\delta^c_n=\min_{K\in \mbox{$\mathcal{K}^n$}}\delta^c(K),\qquad \delta_n^t=
\min_{K\in \mbox{$\mathcal{K}^n$}}\delta^t (K),
\qquad\delta_n^l=\min_{K\in \mbox{$\mathcal{K}^n$}}\delta^l(K),$$

$$\delta^{c\bullet  }_n=\min_{C\in \mbox{$\mathcal{C}^n$}}\delta^c(C),\qquad \delta_n^{t\bullet }=
\min_{C\in \mbox{$\mathcal{C}^n$}}\delta^t (C),
\qquad\delta_n^{l\bullet }=\min_{C\in \mbox{$\mathcal{C}^n$}}\delta^l(C).$$

There are many important and interesting problems about $\delta^c(K)$, $\delta^t(K)$ and $\delta^l(K)$
(see \cite{bras05}). For example,

\medskip
\noindent {\bf Problem 1.} Is there an $n$-dimensional convex body
$K$ satisfying $\delta^t(K)\not= \delta^l(K)$ ?

\medskip
\noindent {\bf Problem 2.} Determine the values of $\delta_n^c$, $\delta^t_n$, $\delta^l_n$,
$\delta_n^{c\bullet  }$, $\delta_n^{t\bullet }$ and $\delta_n^{l\bullet }$, and the corresponding extreme
convex bodies.

\medskip
\noindent {\bf Problem 3.} Determine all the convex polytopes $P$ which
satisfying $\delta^c(P)=1$, $\delta^t(P)=1$ or $\delta^l(P)=1$.

\medskip
In the plane, it was proved by Rogers \cite{roge51} in 1951 that $$\delta^t(K)= \delta^l(K)$$ holds for all convex domains, and by F\'ary \cite{fary50} in
1950 that $$\delta^t_2=\delta^l_2={2\over 3}.$$ However, to determine the value of
$\delta_2^{l\bullet }$ turns out to be extremely challenging. It was conjectured by
Reinhardt \cite{rein34} in 1934, as well as by Mahler \cite{mahl47} in
1947, that
$$\delta_2^{t\bullet } =\delta_2^{l\bullet } ={{8-4\sqrt{2}-\ln 2}\over {2\sqrt{2}-1}}$$
with some smooth octagons as extreme domains. Up to now, this
conjecture is still open. When $n\ge 3$, both Problem 1 and Problem
2 are open. As for Problem 3, we only know partial answers for $n\le 4$
(see Schulte \cite{schu93}).

\medskip
In certain sense, covering can be regarded as a dual or a
counterpart of packing. However, as a research subject, covering is
much younger than packing. By the definitions of $\theta^c(K)$, $\theta^t(K)$ and
$\theta^l(K)$, it is easy to see that
$$1\le \theta^c(K)\le \theta^t(K)\le \theta^l(K)$$
holds for all convex bodies. Like the packing case, we define
$$\theta_n^c=\min_{K\in \mbox{$\mathcal{K}^n$}}\theta^c(K),\qquad \theta_n^t=
\min_{K\in \mbox{$\mathcal{K}^n$}}\theta^t (K),
\qquad\theta_n^l=\min_{K\in \mbox{$\mathcal{K}^n$}}\theta^l(K),$$

$$\theta_n^{c\bullet  }=\min_{C\in \mbox{$\mathcal{C}^n$}}\theta^c(C),\qquad \theta_n^{t\bullet }=
\min_{C\in \mbox{$\mathcal{C}^n$}}\theta^t(C),
\qquad\theta_n^{l\bullet }=\min_{C\in \mbox{$\mathcal{C}^n$}}\theta^l(C).$$
Similar to Problem 1 and Problem 2, we have the following basic problems for covering.

\medskip
\noindent {\bf Problem 4.} Is there an $n$-dimensional convex body
$K$ satisfying $\theta^t(K)\not= \theta^l(K)$ ?

\medskip
\noindent {\bf Problem 5.} Determine the values of $\theta_n^c$, $\theta^t_n$, $\theta^l_n$,
$\theta_n^{c\bullet  }$, $\theta^{t\bullet }_n$ and $\theta_n^{l\bullet }$, and the corresponding extreme
convex bodies.

\medskip
Note that the covering analogue of Problem 3 is itself. In 1950, L.
Fejes T\'oth \cite{feje50} proved that
$$\theta^t(C)=\theta^l(C)$$
holds for every two-dimensional centrally symmetric convex domain
$C$. Unfortunately, up to now, this has neither been generalized to
arbitrary two-dimensional convex domains, nor to higher dimensions. As for the known results
for Problem 5, it follows from results of Sas \cite{sas39} and Fejes
T\'oth \cite{feje50} that
$$\theta_2^{t\bullet }=\theta_2^{l\bullet }={{2\pi }\over {\sqrt{27}}}$$
and ellipses are the only extreme domains; it follows from results of
F\'ary \cite{fary50} and Januszewski \cite{janu10} that
$$\theta_2^t=\theta_2^l = {3\over 2}$$
and triangles are the only extreme domains. Like Problems 1 and 2,
both Problems 4 and 5 are open for $n\ge 3$.

There are many results on $\delta^t(B^n)$, $\delta^l (B^n)$, $\theta^t(B^n)$ and  $\theta^l (B^n)$, on the bounds of $\delta_n^c$, $\delta_n^t$, $\delta_n^l$, $\delta_n^{c\bullet  }$, $\delta^{t\bullet }_n$, $\delta^{l\bullet }_n$, $\theta_n^c$, $\theta_n^t$, $\theta_n^l$, $\theta^{c\bullet  }_n$, $\theta^{t\bullet }_n$ and $\theta_n^{l\bullet }$. Since they are not much
relevant to our purpose, we will not review them here.

Packing, covering and tiling is a research area of mathematics that rich in
challenging problems and fascinating results. For example, the problems and
results about packing densities, covering densities, kissing numbers,
Hadwiger's covering numbers, Rogers' packing-covering constants and
Borsuk's partition numbers. The goal of this paper is to create
a theoretic structure which can be applied to all these problems. Namely, we
will study the geometric structures of $\mathcal{K}^n$ and
$\mathcal{C}^n$ for some particular metrics, and then treat
these numbers as functionals defined on these spaces. We will study
the supderivatives and the integral sums of these functionals over
the metric spaces.

\bigskip
\section*{2. Spaces of Convex Bodies}

\medskip
Let $K_1+K_2$ denote the {\it Minkowski sum} of $K_1$ and $K_2$
defined by
$$K_1+K_2=\{ {\bf x}_1+{\bf x}_2:\ {\bf x}_i\in K_i\},$$
let $\| \cdot \|^*$ denote the {\it Hausdorff metric} on
$\mathcal{K}^n$ defined by
$$\| K_1, K_2\|^*=\min \left\{ r:\ K_1\subseteq K_2+rB^n,\ K_2\subseteq
K_1+rB^n\right\},$$ and let $\{ \mathcal{K}^n, \| \cdot \|^*\}$
denote the space of $\mathcal{K}^n$ with metric $\| \cdot \|^*$.
It is easy to see that, for $\lambda_i\in \mathbb{R}$ and $K_i\in
\mathcal{K}^n$,
$$\lambda_1K_1+\lambda_2K_2+\ldots +\lambda_mK_m\in \mathcal{K}^n.$$
In certain sense, the space $\mathcal{K}^n$ has linear structure.

\medskip
In 1916, Blaschke proved the following theorem:

\smallskip \noindent {\bf Blaschke's selection theorem.} {\it Let
$r_1$ and $r_2$ be two positive numbers with $r_1<r_2$. For
any infinite sequence of $n$-dimensional convex bodies $\{
K_1, K_2, K_3, \ldots \}$ all satisfying $r_1B^n\subseteq
K_i\subseteq r_2B^n$, there is a subsequence $\{ K'_1, K'_2, K'_3,
\dots \}$ and an $n$-dimensional convex body $K_0$ satisfying}
$$\lim_{i\to\infty}\| K'_i, K_0\|^*=0.$$

This theorem guarantees the local compactness of $\{ \mathcal{K}^n,
\|\cdot \|^*\}$. It is easy to show that all $\delta^c(K)$, $\delta^t(K)$,
$\delta^l(K)$, $\theta^c(K)$, $\theta^t(K)$ and $\theta^l(K)$ are bounded continuous
functionals defined on $\{ \mathcal{K}^n, \|\cdot \|^*\}$. However, the Hausdorff metric
has a disadvantage that it can not distinguish the shapes of the
convex bodies. Let $r$ be a positive number. Clearly,
$K$ and $rK$ have the same shape, and
$$\delta^t(rK)=\delta^t(K),\qquad \delta^l (rK)=\delta^l (K),\qquad
\theta^t(rK)=\theta^t(K), \qquad \quad \theta^l (rK)=\theta^l (K),$$
$$\tau^t (rK)=\tau (K),\qquad \tau^l (rK)=\tau^* (K),\qquad
h(rK)=h(K) \quad\  and \quad\ b(rK)=b(K).$$
However, on the other hand, it can be easily shown that $\| K, rK\|^*$ can
be arbitrary large when $r\to \infty$. This shows the disadvantage of the Hausdorff metric
in the study of these numbers. There are several other metrics (see \cite{grub93}),
one of them is particular important for our purpose: the {\it Banach-Mazur metric}.

\medskip
In 1948, generalizing a two-dimensional result of Behrend, John \cite{john48}
proved the following basic theorem:

\smallskip\noindent {\bf John's Theorem.} {\it For each
$n$-dimensional convex body $K$ there is an ellipsoid $E$ satisfying
$E\subseteq K\subseteq nE$; For each $n$-dimensional centrally
symmetric convex body $C$ there is an ellipsoid $E'$ satisfying
$E'\subseteq C\subseteq \sqrt{n}E'.$}

\medskip
This theorem sparked the idea of reduction. Let $\mathcal{T}^n$ denote the family of all nonsingular affine linear transformations from $\mathbb{E}^n$ to $\mathbb{E}^n$, and let $\| \cdot \|$ denote
the Banach-Mazur metric defined by
$$\| K_1, K_2\| =\log\ \min \left\{ r: \ K_1\subseteq \sigma (K_2)\subseteq
rK_1+{\bf x}; \ {\bf x}\in \mathbb{E}^n; \ \sigma\in \mathcal{T}^n\right\}.$$ It
is known (easy to prove) that both $\{ \mathcal{K}^n, \| \cdot \|\}$
and $\{ \mathcal{C}^n, \|\cdot \|\}$ are metric spaces.  Let $\|
X\|$ denote the diameter of a set $X$ with respect to the
Banach-Mazur metric. It follows by John's theorem and the triangular
inequality of $\|\cdot \|$ that
$$\| \mathcal{K}^n\| \le 2\log n\qquad and \qquad  \| \mathcal{C}^n\| \le \log
n.$$ Therefore, by John's theorem and Blaschke's selection theorem,
both $\{\mathcal{K}^n, \|\cdot \|\}$ and $\{\mathcal{C}^n, \|\cdot \|\}$ are bounded, connected and
compact. This is essentially different from the Hausdorff metric.
Since each centrally symmetric convex body corresponds to a Banach
space, the following problem (see \cite{tomc89}) is fundamental in
Functional Analysis, as well as in Convex Geometry.

\medskip
\noindent {\bf Problem 6.} Determine the values of $\|
\mathcal{K}^n\| $ and $\| \mathcal{C}^n\|$.

\medskip
Let $I^2$ denote a square, let $H$ denote a regular hexagon and define
$$D=\left\{ (x,y):\ |y|\le 1,\
x^2+y^2\le 2,\ \mbox{${1\over 2}$}x^2+y^2\le \mbox{${4\over 3}$}\right\}.$$ In 1981,
Stromquist \cite{stro81} proved that for all $C\in \mathcal{C}^2$ we
have
$$\| C, D\|\le \mbox{$1\over 2$}(\log
3-\log 2),$$ where equality holds if and only if $C=I^2$ or $C=H$.
Therefore, combined with the fact that $$\| I^2, H\|=\log 3-\log 2\eqno (1)$$
which was discovered by Asplund \cite{aspl60} in 1960, we get
$$\| \mathcal{C}^2\|=\log 3-\log 2.\eqno (2)$$ Up to now, this is the only known exact answer to Problem 6.

In $\{\mathcal{K}^n, \| \cdot
\|\}$, it can be shown that
$$\| K, K'\|=0$$
if and only if $K'=\sigma(K)$ for some $\sigma\in\mathcal{T}^n$. This
observation leads to another representation of $\{ \mathcal{K}^n,
\|\cdot \|\}$. Namely, there is a bounded, connected and compact subset
(even in the sense of the Hausdorff metric)
$\widehat{\mathcal{K}^n}$ of $\mathcal{K}^n$ such that
$$\mathcal{K}^n=\widehat{\mathcal{K}^n}\otimes \mathcal{T}^n.$$
Similarly, there is a bounded, connected and compact subset
$\widehat{\mathcal{C}^n}$ of $\mathcal{C}^n$ such that
$$\mathcal{C}^n=\widehat{\mathcal{C}^n}\otimes \mathcal{T}^n.$$
Therefore, the relation between $\{\mathcal{K}^n, \|\cdot \|^*\}$ and $\{\mathcal{K}^n, \|\cdot\|\}$,
as well as $\{\mathcal{C}^n, \|\cdot \|^*\}$ and $\{\mathcal{C}^n, \|\cdot\|\}$, is similar to the relation between $\mathbb{E}^n$ and the spherical space $\partial (B^n)$ with the spherical metric.

\medskip\noindent
{\bf Remark 1.} Usually, a metric $d(\cdot )$ requires that $d({\bf x}, {\bf y})=0$
if and only if ${\bf x}={\bf y}$. However, this is not true in the
cases of $\{ \mathcal{K}^n, \| \cdot \|\}$ and $\{ \mathcal{C}^n, \|
\cdot \|\}$. In these spaces, $K\in \mathcal{F}$ means
$$\{ K':\ K'=\sigma(K), \ \sigma\in \mathcal{T}^n\} \subseteq \mathcal{F}.$$

Let $\rho $ be a small positive number and let $K$ be an $n$-dimensional convex body. We call
$$\mathcal{B}(K,\rho )=\left\{ K':\ K'\in\mathcal{K}^n,\ \|K, K'\|\le \rho \right\}$$ a
ball in $\{\mathcal{K}^n, \|\cdot \|\}$ centered at $K$ of radius
$\rho$. Just like the Euclidean case, we call $\{ K':\
K'\in\mathcal{K}^n,\ \|K, K'\|= \rho \}$ and $\{ K':\
K'\in\mathcal{K}^n,\ \|K, K'\|< \rho \}$ the boundary and the
interior of $\mathcal{B}(K,\rho )$, and denote them by $\partial
(\mathcal{B}(K,\rho ) )$ and ${\rm int} (\mathcal{B}(K,\rho ))$,
respectively. Open sets in $\{\mathcal{K}^n, \|\cdot \|\}$ can be defined in a routine way.
Similar concepts can be defined in $\{\mathcal{C}^n, \|\cdot \|\}$, $\{\mathcal{K}^n, \|\cdot \|^*\}$,
$\{\mathcal{C}^n, \|\cdot \|^*\}$ and etc.

\medskip
Now, we are facing the following fundamental problem:

\medskip
{\it Are there geometrical useful measures
on $\{ \mathcal{K}^n, \|\cdot \|^*\}$, $\{ \mathcal{C}^n, \|\cdot \|^*\}$, $\{ \mathcal{K}^n,
\|\cdot \|\}$ and $\{ \mathcal{C}^n, \|\cdot \|\}$?}

\medskip
If the answer is \lq\lq yes", it would provide powerful tools to study the functionals $\delta^c(K)$, $\delta^t(K)$, $\delta^l(K)$, $\theta^c(K)$, $\theta^t(K)$, $\theta^l(K)$, $\tau^c(K)$, $\tau^t(K)$, $\tau^l(K)$ and etc defined on these spaces. Unfortunately, in 1986 Bandt and Baraki \cite{band86} proved the following result, which in certain sense gave a negative answer to this problem.  {\it When $n\ge 2$, there is no positive $\delta$-finite Borel measure on $\{ \mathcal{K}^n, \|\cdot \|^*\}$ which is invariant with respect to all isometries in it.}

\medskip\noindent
{\bf Remark 2.} It was proved by Gruber and Lettl \cite{grub80} that, $\sigma $ is an isometry in $\{ \mathcal{K}^n,
\|\cdot \|^*\}$ if and only if
$$\sigma (K)=\varsigma (K)+K_0,$$
where $\varsigma $ is a rigid motion in $\mathbb{E}^n$, $K_0$ is a convex body and $+$ is the Minkowski sum.

\medskip
In 2010, Hoffmann \cite{hoff10} constructed the following measure on $\{ \mathcal{K}^n,
\|\cdot \|^*\}$: Let $(K_i)_{i\in \mathbb{N}}$ be a sequence of convex bodies which is dense in $\{ \mathcal{K}^n,
\|\cdot \|^*\}$, let $(\alpha_i)_{i\in \mathbb{N}}$ be a sequence of positive number such that $\sum \alpha_i <\infty $, and for $K\in \mathcal{K}^n$, let $\delta_K$ denote the {\it Dirac measure} concentrated at $K$. For any open set
$\mathcal{O}$ of $\{ \mathcal{K}^n, \|\cdot \|^*\}$ we define its measure by
$$\mu (\mathcal{O})=\sum_{K_j\in \mathcal{O}}\alpha_j\delta_{K_j}.$$

In $\{ \mathcal{C}^n, \|\cdot \|^*\}$, $\{ \mathcal{K}^n, \|\cdot \|\}$ and $\{ \mathcal{C}^n, \|\cdot \|\}$
one can do the similar constructions as well. However, such measures seem not geometrically useful.

\medskip\noindent
{\bf Definition 1.} {\it Let $\beta$ be a positive number. A subset $\mathcal{X}$ of
$\mathcal{K}^n$ is called a $\beta$-net in $\{ \mathcal{K}^n,
\|\cdot \|\}$ if for each $K\in \mathcal{K}^n$ there is a $K'\in \mathcal{X}$
satisfying $$\| K, K'\|<\beta .$$ We define $\ell(n,\beta )$ to be the
smallest number of convex bodies which forms a $\beta $-net in $\{ \mathcal{K}^n,
\|\cdot \|\}$.}

\medskip
In fact, $\ell(n,\beta )$ is the smallest number of
open balls of radius $\beta $ that their union can cover the
whole space $\{ \mathcal{K}^n, \|\cdot \|\}$. In 2010 Zong \cite{zong10}
proved the following result: {\it The minimal cardinality of
$\beta$-nets in $\{\mathcal{K}^n,\| \cdot \|\}$ is bounded by
$$\ell (n, \beta )\le
\left\lfloor\mbox{${{7n}\over \beta }$}\right\rfloor^{c\cdot
14^n\cdot n^{2n+3}\cdot \beta^{-n}},\eqno (3)$$ where $c$ is a suitable
constant.}

\medskip
This bound is far from sharp. Nevertheless, it reveals the fact that
$\ell(n,\beta )$ is bounded from above.

\medskip\noindent
{\bf Definition 2.} {\it We say a family of balls
$\mathcal{F}=\{\mathcal{B}(K_i,\rho ):\ i=1, 2, \ldots, m\}$ forms a
packing in $\{\mathcal{K}^n, \|\cdot \|\}$ if
$${\rm int} (\mathcal{B}(K_i,\rho ))\cap {\rm int} (\mathcal{B}(K_j,\rho
))=\emptyset, \quad i\not= j,$$ and define $m(n,\rho )$ to be the
maximal number of balls of radius $\rho $ which can be packed into
$\{\mathcal{K}^n, \|\cdot \|\}$.}

\medskip
Recall that $\ell(n,\beta )$ is the minimal number of balls of radius
$\beta$ that the union of their interiors covers $\{\mathcal{K}^n,
\|\cdot \|\}$. It can be deduced (see \cite{lore66}) that, for any
positive number $\omega$,
$$m(n,\omega )\le \ell(n,\omega )\le m\left(n,\mbox{$1\over 2$}\omega \right).\eqno(4)$$
The first inequality can be deduced from the fact that, if $\{
\mathcal{B}(K_i,\omega ):\ i=1, 2, \ldots, m\}$ forms a ball packing in
$\{ \mathcal{K}^n, \| \cdot \| \}$ and $\{ {\rm
int}(\mathcal{B}(K'_j,\omega )):\ j=1, 2, \ldots, l\}$ forms a covering
of $\{ \mathcal{K}^n, \| \cdot \| \}$, then each ${\rm
int}(\mathcal{B}(K'_j, \omega ))$ contains at most one $K_i$ and
therefore $l\ge m$. Otherwise, if ${\rm int} (\mathcal{B}(K'_j,\omega
))$ contains both $K_{i_1}$ and $K_{i_2}$, then we have $$K'_j\in
{\rm int}(\mathcal{B}(K_{i_1}, \omega ))\cap {\rm
int}(\mathcal{B}(K_{i_2},\omega )),$$ which contradicts the assumption
that $\{ \mathcal{B}(K_i, \omega ):\ i=1, 2, \ldots , m\}$ forms a
packing in $\{ \mathcal{K}^n, \|\cdot\|\}.$

The second inequality in (4) can be shown by the fact that, if $\{
\mathcal{B}(K_i,{1\over 2}\omega ):\ i=1, 2, \ldots, m'\}$ forms a
packing of maximal cardinality in $\{ \mathcal{K}^n, \| \cdot \|
\}$, then $\{ {\rm int}(\mathcal{B}(K_i,\omega )):\ i=1, 2, \ldots,
m'\}$ forms a covering of $\{ \mathcal{K}^n, \|\cdot\|\}$ and
therefore $m'\ge \ell(n,\omega )$. Otherwise, if
$$K\not\in \bigcup_{i=1}^{m'}{\rm int} (\mathcal{B}(K_i, \omega )),$$
then we have
$${\rm int}(\mathcal{B}(K, \mbox{${1\over 2}$}\omega ))\cap
{\rm int}(\mathcal{B}(K_i, \mbox{${1\over 2}$}\omega ))=\emptyset ,
\quad i=1,2, \ldots , m',$$ which contradicts the maximum assumption
on $m'$.

It follows by (3) and (4) that $m(n, \rho)$ is bounded from
above as well. The following problem is basic for the structure of
$\{ \mathcal{K}^n, \| \cdot \| \}$. Clearly, similar question can be
asked for $\{ \mathcal{C}^n, \| \cdot \| \}$.

\medskip
\noindent {\bf Problem 7.} For a given dimension $n$ and some particular $\beta <\|\mathcal{K}^n\|$, determine the values of $\ell(n,\beta )$ and $m(n,\beta )$; For a given dimension $n$ and small $\beta $, determine (or
estimate) the asymptotic orders of $\ell(n,\beta )$ and $m(n,\beta )$ when $\beta\to 0$.

\medskip
In $\{ \mathcal{K}^n, \| \cdot \|^* \}$ one can similarly define ball coverings and ball packings. Since
$\|\mathcal{K}^n \|^*=\infty$, one can't define analogues of $\ell (n, \beta )$ and $m(n, \rho )$ in $\{ \mathcal{K}^n, \| \cdot \|^* \}$. However, we can define ball coverings, ball packings and analogues of $\ell (n, \beta )$ and $m(n, \rho )$ in $$\mathcal{K}^{n*}=\{ K\in \mathcal{K}^n:\ B^n\subseteq K\subseteq nB^n \}$$
with the metric $\|\cdot \|^*$. Let $\ell^* (n, \beta )$ and $m^*(n, \rho )$ denote the analogues of $\ell (n, \beta )$ and $m(n, \rho )$ in $\mathcal{K}^{n*}$ with respect to the Hausdorff metric, respectively. Similar to (4), for any positive number $\omega $, we have
$$m^*(n, \omega )\le \ell^*(n, \omega )\le m^*(n, \mbox{$1\over 2$}\omega )\eqno (5)$$
as well.

From the intuitive point of view, it is easy to imagine that $\| K_1, K_2\|/\|K_1, K_2\|^*$ can be arbitrarily small. In fact, it can be arbitrarily large as well. Let $I^2$ be the unit square, let $H$ be the regular hexagon with unit edge, and let $\epsilon $ be a small positive number. Then we have
$$\|\epsilon I^2, \epsilon H\| =\log 3-\log 2,$$
$$\|\epsilon I^2, \epsilon H\|^* \le (1-{{\sqrt{2}}/ 2})\epsilon$$
and
$${{\|\epsilon I^2, \epsilon H\| }\over {\|\epsilon I^2, \epsilon H\|^*}}\ge {{2(\log 3-\log 2)}\over {2-\sqrt{2}}}\cdot
{1\over \epsilon },$$
which can be arbitrarily large when $\epsilon \to 0$. Nevertheless, we have the following result which reflects the relation between $\|\cdot\|$ and $\|\cdot\|^*$.

\medskip
\noindent
{\bf Theorem 1.} {\it For every pair of convex bodies $K_1$ and $K_2$ in $\mathcal{K}^{n*}$, we have}
$$\|K_1, K_2\|\le 2\cdot \| K_1, K_2\|^*.$$

\smallskip\noindent
{\bf Proof.} Assume that $\| K_1, K_2\|^*=r^*$. Then we have
$$K_1\subseteq K_2+r^*B^n$$
and
$$K_2\subseteq K_1+r^*B^n.$$
On the other hand, since both $K_1$ and $K_2$ belong to $\mathcal{K}^{n*}$, we have
$$B^n\subseteq K_i\subseteq nB^n,\qquad i=1,\ 2.$$
Thus, we have
$$K_1\subseteq K_2+r^*K_2=(1+r^*)K_2,$$
$$K_2\subseteq K_1+r^*K_1=(1+r^*)K_1,$$
$$\mbox{${1\over {1+r^*}}$}K_1\subseteq K_2\subseteq (1+r^*)K_1,$$
$$\mbox{${1\over {1+r^*}}$}K_2\subseteq K_1\subseteq (1+r^*)K_2,$$
and
$$\| K_1, K_2\|\le 2\log (1+r^*)\le 2r^*=2\| K_1, K_2\|^*.$$
The theorem is proved. \hfill{$\square$}

\bigskip
\section*{3. Functionals on $\mathcal{K}^n$ and $\mathcal{C}^n$}

\medskip
By routine arguments it can be shown that $\delta^c(K)$, $\delta^t(K)$, $\delta^l(K)$,
$\theta^c(K)$, $\theta^t(K)$ and $\theta^l(K)$ are continuous functionals defined on $\{ \mathcal{K}^n,
\|\cdot \|^*\}$, and $\delta^{c\bullet  }(C)$, $\delta^{t\bullet }(C)$, $\delta^{l\bullet }(C)$,
$\theta^{c\bullet  }(C)$, $\theta^{t\bullet }(C)$ and $\theta^{l\bullet }(C)$ are continuous in $\{ \mathcal{C}^n,
\|\cdot \|^*\}$. Similarly, $\delta^t(K)$, $\delta^l(K)$, $\theta^t(K)$ and $\theta^l(K)$ are continuous in $\{ \mathcal{K}^n, \|\cdot \|\}$, and $\delta^{t\bullet }(C)$, $\delta^{l\bullet }(C)$, $\theta^{t\bullet }(C)$ and $\theta^{l\bullet }(C)$ are continuous in $\{ \mathcal{C}^n, \|\cdot \|\}$. However, $\tau^c(K)$, $\tau^t(K)$, $\tau^l(K)$, $h(K)$ and $b(K)$ are not continuous in $\{ \mathcal{K}^n, \|\cdot \|^*\}$ and $\delta^c(K)$, $\theta^c(K)$, $\tau^c(K)$, $\tau^t(K)$, $\tau^l(K)$, $h(K)$ and $b(K)$ are not continuous in $\{ \mathcal{K}^n, \|\cdot \|\}$.

Just like the real functions defined in $\mathbb{R}$, if $f(K)$ and $g(K)$ are continuous in
$\{ \mathcal{K}^n, \|\cdot \|^*\}$, then both $f(K)+g(K)$ and $f(K)\cdot g(K)$ are continuous. Of course,
analogues are also true in $\{ \mathcal{C}^n, \|\cdot \|^*\}$, $\{ \mathcal{K}^n,
\|\cdot \|\}$, $\{ \mathcal{C}^n, \|\cdot \|\}$ and similar metric spaces.

Let $\sigma$ denote a non-singular affine linear transformation from $\mathbb{E}^n$ to
$\mathbb{E}^n$.  It is easy to see that $\sigma (K)$ is a convex body provided $K$
is such one, and
$$\delta^t(\sigma (K))=\delta^t(K),\qquad \delta^l (\sigma (K))=\delta^l (K),\qquad
\theta^t(\sigma (K))=\theta^t(K),$$
$$\theta^l (\sigma (K))=\theta^l (K), \qquad \tau^t(\sigma (K))=\tau^t (K),\qquad \tau^l(\sigma (K))=\tau^l (K).$$
In other words, we have
$$\delta^t(K_1)=\delta^t(K_2),\qquad\delta^l(K_1)=\delta^l(K_2), \qquad\theta^t(K_1)=\theta^t(K_2),$$
$$\theta^l(K_1)=\theta^l(K_2),\qquad \tau^t(K_1)=\tau^t(K_2),\qquad \tau^l(K_1)=\tau^l(K_2)$$ whenever $\| K_1, K_2\|=0$.
Therefore, we can treat $\delta^t(K)$, $\delta^l(K)$, $\theta^t(K)$, $\theta^l(K)$, $\tau^t(K)$ and $\tau^l(K)$ as functionals defined on $\{\mathcal{K}^n, \|\cdot \| \}$.

Let $T_3$ denote a regular tetrahedron and let $I^3$ denote a unit cube in $\mathbb{E}^3$. Let ${\bf v}_1$
and ${\bf v}_2$ be two opposite vertices of $I^3$ and enumerate the other vertices of $I^3$ as ${\bf v}_3$,
${\bf v}_4$, $\ldots$, ${\bf v}_8$ such that ${\bf v}_i{\bf v}_{i+1}$ are edges of $I^3$. Of course, here ${\bf v}_9={\bf v}_3$. Then $I^3$ can be triangulated into six congruent tetrahedra ${\bf v}_1{\bf v}_2{\bf v}_i{\bf v}_{i+1}$, $i=3,$ $4,$ $\ldots ,$ $8$. For convenience, let $T'_3$ denote the tetrahedron ${\bf v}_1{\bf v}_2{\bf v}_3{\bf v}_4$ and let $\sigma $ be an affine linear transformation such that $\sigma (T_3)=T'_3$. It is known that
$\mathbb{E}^3$ can be tiled by $I^3$, and therefore also by $T'_3$, but can't be tiled by $T_3$. In other words,
we have $\delta^c(T_3)<1$, $\theta^c(T_3)>1$ and $\delta^c(T'_3)=\theta^c(T'_3)=1$. Consequently, we have (see \cite{laga13})
$$\delta^c(\sigma (T_3))\not= \delta^c(T_3)$$
and
$$\theta^c(\sigma (T_3))\not= \theta^c(T_3).$$
Therefore, to study $\delta^c(K)$ and $\theta^c(K)$, we have to work in $\{\mathcal{K}^n, \|\cdot \|^*\}$.

\medskip
There are several approaches to study the relations between these functionals. For example, in 1950, Rogers \cite{roge50} proved that
$$\theta^l(C)\le 3^n\cdot\delta^l(C)\eqno (6)$$
and
$$\theta^t(C)\le 2^n\cdot\delta^t(C)\eqno (7)$$
hold for all $C\in \mathcal{C}^n$. To this end, he introduced and studied the lattice packing-covering constant $\phi^l (C)$. For a lattice $\Lambda $, let $\rho (C, \Lambda )$ denote the largest number $\rho$
such that $\rho C+\Lambda$ is a packing and let $\rho'(C,\Lambda )$ denote the smallest number $\rho'$ such that $\rho'C+\Lambda $ is a covering of $\mathbb{E}^n$. Then we define
$$\phi^l (C)=\min_\Lambda {{\rho'(C,\Lambda )}\over {\rho (C,\Lambda )}}$$
and call it the lattice packing-covering constant of $C$.
Similarly, one can define the translative packing-covering constant $\phi^t(C)$. Clearly both $\phi^l(C)$ and $\phi^t(C)$ are affinely invariant continuous functionals defined on $\{\mathcal{C}^n, \|\cdot \|\}$.
In fact, (6) and (7) can be deduced from
$$\theta^l(C)\le \phi^l(C)^n\cdot\delta^l(C)\eqno (8)$$
and
$$\theta^t(C)\le \phi^t(C)^n\cdot\delta^t(C),\eqno (9)$$
respectively. In $\mathbb{E}^2$ it was shown by Zong \cite{zong08} that
$$\phi^t(C)=\phi^l(C)\le 2(2-\sqrt{2}),$$
where the equality holds if and only if $C$ is an affine regular octagon.

In 2001, Ismailescu \cite{isma01} proved that
$$1-\delta^l(K)\le \theta^l(K)\le 1.25\sqrt{1-\delta^l(K)}\eqno(10)$$
holds for all $K\in \mathcal{K}^2$. Results such as (8), (9) and (10) can be regarded as examples to study relations of particular functionals defined in $\mathcal{C}^n$ and $\mathcal{K}^n$, respectively.

\medskip
Both Hadwiger's covering number $h(K)$ and Borsuk's partition number $b(K)$ are discontinuous in $\{ \mathcal{K}^n, \|\cdot\|^*\}$. Let $m$ be a fixed positive integer. In 2010, Zong \cite{zong10} introduced and studied two functionals
$\gamma_m(K)$ and $\varphi_m(K)$. Namely, $\gamma_m(K)$ is the smallest number $r$ such that $K$ can be covered by
$m$ translates of $rK$ and $\varphi_m(K)$ is the smallest number $\mu$ such that $K$ can be divided into $m$ subsets
$X_1$, $X_2$, $\cdots$, $X_m$ such that
$$d(X_i)\le \mu\cdot d(K)$$
holds for all the subsets, where $d(X)$ denote the Euclidean diameter of $X$. It was proved that Hadwiger's conjecture is equivalent with
$$\gamma_{2^n}(K)\le c_1<1$$
holds for all $K\in \mathcal{K}^n$, where $c_1$ is a suitable positive constant, and Borsuk's conjecture is equivalent with
$$\varphi_{n+1}(K)\le c_2<1$$
holds for all $K\in \mathcal{K}^n$, where $c_2$ is a suitable positive constant. It is important that both $\gamma_m(K)$ and $\varphi_m(K)$ are continuous in $\{\mathcal{K}^n, \|\cdot \|^*\}$.

\medskip
Next, let us make a couple of observations which show some importance of studying these functionals.

\medskip\noindent
{\bf Observation 1.} As shown by (2) that the diameter of $\{ \mathcal{C}^2, \|\cdot\|\}$ is $\log 3-\log 2$. In $\mathbb{E}^2$ it is well-known that
$$\delta^t(C)=\delta^l(C)=\theta^t(C)=\theta^l(C)=1\eqno (11)$$
if and only if $C$ is a centrally symmetric hexagon or a centrally symmetric parallelogram. Let $\mathcal{X}$ denote this set. It is interesting to note that $\mathcal{X}$ is a connect compact subset of $\{ \mathcal{C}^2, \| \cdot\|\}$ without interior point, and by (1) and (2)
$$\| \mathcal{X} \|= \| \mathcal{C}^2 \|= \log 3-\log 2.$$
On the other hand, both $\theta^t(C)$ and $\theta^l(C)$ attain their maxima if and only if $C$ is an ellipse, a zero diameter set in $\{\mathcal{C}^2, \|\cdot\|\}$. To determine the minima of $\delta^t(C)$ and $\delta^l(C)$ and the corresponding extreme domains is still a challenging open problem. Similar observation can be made for $\delta^c(K)$ and $\theta^c(K)$ in $\{ \mathcal{K}^2, \| \cdot\|^*\}$. But the situation can be much more complicated.

\medskip\noindent
{\bf Observation 2.} It was proved by Minkowski that
$$\tau^l(K)\le 3^n-1,$$
where the upper bound can be attained if and only if $K$ is an $n$-dimensional parallelopiped, and
$$\tau^l(K)\le 2(2^n-1)$$
whenever $K$ is an $n$-dimensional strict convex body. As a consequence, since the set of strict convex bodies is dense in $\{\mathcal{K}^n, \|\cdot\|\}$, the set $\{ K\in \mathcal{K}^n:\ 2(2^n-1)<\tau^l(K)\le 3^n-1\}$ contains no open subset in $\{ \mathcal{K}^n, \|\cdot \|\}$.

\medskip
Now, we reformulate and generalize Problems 1 and 2 into the following functional form.

\medskip
\noindent
{\bf Problem 8.} To study $\delta^c(K)-\delta^t(K)$, $\delta^c(K)-\delta^l(K)$, $\delta^t(K)-\delta^l(K)$, $\theta^c(K)-\theta^t(K)$, $\theta^c(K)-\theta^l(K)$, $\theta^t(K)-\theta^l(K)$,
$\tau^t(K)-\tau^l(K)$, $\theta^t(K)/\delta^t(K)$, $\theta^l(K)/\delta^l(K)$ and other similar functionals in $\{\mathcal{K}^n, \|\cdot\|^*\}$ or in $\{\mathcal{K}^n, \|\cdot\|\}$.

\bigskip
\section*{4. Supderivatives}
\medskip
\noindent
{\bf Definition 3.} {\it Let $f(K)$ be a continuous functional defined on $\{\mathcal{K}^n, \|\cdot\|\}$
and let $K_0$ be a particular $n$-dimensional convex body. If
$$f'(K_0): =\limsup_{K\to K_0}{{|f(K)-f(K_0)|}\over {\| K, K_0\|}}$$
is finite, we call it the supderivative of $f(K)$ at $K_0$ with respect to $\|\cdot \|$. Similarly, one can define the supderivative $f^*(K_0)$ of $f(K)$ at $K_0$ with respect to $\| \cdot \|^*$.}

\medskip
Supderivatives, like the derivatives, have some basic properties. For example, it can be shown that
$$\left( f(K)+g(K)\right)'\le f'(K)+g'(K),$$
$$\left( f(K)\cdot g(K)\right)'\le f'(K)\cdot |g(K)|+ |f(K)|\cdot g'(K),$$
and
$$\left({{f(K)}\over {g(K)}}\right)'\le {{f'(K)\cdot |g(K)|+ |f(K)|\cdot g'(K)}\over {g(K)^2}}.$$

\medskip
Now, we estimate the supderivatives of some particular functionals.

\medskip\noindent
{\bf Theorem 2.} {\it Let $f(K)$ to be $\theta^c(K)$, $\theta^t(K)$, $\theta^l(K)$, $\delta^c(K)$, $\delta^t(K)$ or $\delta^l(K)$, and assume that $K_0\in \mathcal{K}^{n*}$. Then, for $K\in \mathcal{K}^{n*}$, we have}
$$f(K)\le (1+\| K, K_0\|^*)^{2n}f(K_0).$$

\medskip\noindent
{\bf Proof.} For convenience, we write
$$r=\| K, K_0\|^*.$$
Then, by the definitions of $\|\cdot\|^*$ and $\mathcal{K}^{n*}$ we have
$$K_0\subseteq K+rB^n,$$
$$B^n\subseteq K$$
and therefore
$$K_0\subseteq K+rK=(1+r)K.\eqno(12)$$
Similarly, we also have
$$K\subseteq (1+r)K_0.\eqno(13)$$

Let $\theta (\Im )$ denote the density of a covering system $\Im $ of $\mathbb{E}^n$. Assume that
$$\Im_0=\{K_0+{\bf x}_i:\ i=1, 2, \ldots \}$$
is a translative covering of $\mathbb{E}^n$ with the density
$$\theta (\Im_0)=\theta^t(K_0).$$
Then, by (12) it is easy to see that
$$\Im_1=\{ (1+r)K+{\bf x}_i:\ i=1, 2, \ldots \}$$
is a translative covering of $\mathbb{E}^n$ with density
$$\theta (\Im_1)\ge\theta^t(K).\eqno (14)$$
Furthermore, by (13) and (14), the system
$$\Im_2=\{(1+r)^2K_0+{\bf x}_i:\ i=1, 2, \ldots \}$$
covers $\mathbb{E}^n$ with density $\theta (\Im_2 )$, which satisfies
$$\theta (\Im_2)=(1+r)^{2n}\theta^t(K_0)\ge \theta (\Im_1)\ge \theta^t(K).$$
The $f(K)=\theta^t(K)$ case of the theorem is proved. The other covering cases can be dealt with by similar arguments.

Let $\delta (\wp )$ denote the density of a packing system $\wp $ in $\mathbb{E}^n$. Assume that
$$\wp_0=\{ (1+r)^2K+{\bf x}_i:\ i=1, 2, \ldots \}$$
is a translative packing in $\mathbb{E}^n$ with the density
$$\delta (\wp_0)=\delta^t(K).$$
Then, by (12) it is easy to see that
$$\wp_1=\{ (1+r)K_0+{\bf x}_i:\ i=1, 2, \ldots \}$$
is a translative packing in $\mathbb{E}^n$ with density
$$\delta (\wp_1)\le \delta^t(K_0).\eqno (15)$$
Furthermore, by (13) and (15), the system
$$\wp_2=\{K+{\bf x}_i:\ i=1, 2, \ldots \}$$
is a packing in $\mathbb{E}^n$ with density $\delta (\wp_2 )$, which satisfies
$$\delta (\wp_2)=(1+r)^{-2n}\delta^t(K)\le \delta (\wp_1)\le \delta^t(K_0).$$
Thus, we get
$$\delta^t(K)\le (1+r)^{2n}\delta^t(K_0),$$
which proves the $f(K)=\delta^t(K)$ case of the theorem. The other packing cases can be dealt with by similar arguments. \hfill{$\square$}

\medskip\noindent
{\bf Corollary 1.} {\it Let $f(K)$ to be $\theta^c(K)$, $\theta^t(K)$,  $\theta^l(K)$, $\delta^c(K)$, $\delta^t(K)$ or $\delta^l(K)$. There are two positive constants $c^*(n)$ and $d^*(n)$, which depend on $n$, such that
$$|f(K_1)-f(K_2)|\le c^*(n)\cdot \| K_1, K_2\|^*$$
holds for any pair $K_1$, $K_2\in \mathcal{K}^{n*}$, provided $\| K_1, K_2\|^*\le d^*(n)$. In particular for all $K\in \mathcal{K}^{n*}$, we have}
$$f^*(K)\le 2n.$$

\medskip\noindent
{\bf Remark 3.} When $n$ is small, both $c^*(n)$ and $d^*(n)$ can be taken precisely. In 1950 it was shown by F\'ary \cite{fary50} that
$$\theta^l(K)\le {3\over 2}$$
holds for all two-dimensional convex domains $K$. At the same time, Macbeath \cite{macb51} discovered that, for $K\in \mathcal{K}^n$ there is a cylinder $H$ inscribed in $K$ with
$${\rm vol}(H)\ge {{(n-1)^{n-1}}\over {n^n}}{\rm vol}(K).$$
Consequently, in $\mathbb{E}^3$ one can deduce
$$\theta^l(K)\le {3\over 2}\cdot {{3^3}\over {2^2}}={{3^4}\over {2^3}}.$$
Thus, in the covering cases, one can take $c^*(2)=12$ and $d^*(2)={1\over 3}$, and take $c^*(3)=122$ and $d^*(3)={1\over 4}$; in the packing cases, one can take $c^*(2)=8$ and $d^*(2)={1\over 3}$, and take $c^*(3)=12$ and $d^*(3)={1\over 4}$.

\bigskip\noindent
{\bf Theorem 3.} {\it Let $f(K)$ to be $\theta^t(K)$, $\theta^l(K)$, $\delta^t(K)$ or $\delta^l(K)$, and let $K_0$ be an $n$-dimensional convex body. Then, for any $K\in \mathcal{K}^n$, we have}
$$f(K)\le e^{n\cdot\| K, K_0\| }f(K_0).$$

\medskip\noindent
{\bf Proof.} For convenience, we write
$$r=e^{\| K, K_0\|}.$$
By the definition of $\|\cdot\|$ we have
$$K_0\subseteq \sigma (K)\subseteq rK_0+{\bf x},\eqno(16)$$
where $\sigma $ is a non-singular affine linear transformation from $\mathbb{E}^n$ to $\mathbb{E}^n$
and ${\bf x}$ is a suitable point in $\mathbb{E}^n$.

Assume that
$$\Im_0=\{K_0+{\bf x}_i:\ i=1, 2, \ldots \}$$
is a translative covering of $\mathbb{E}^n$ with the density
$$\theta (\Im_0)=\theta^t(K_0).$$
Then, by (16) it is easy to see that
$$\Im_1=\{ \sigma (K)+{\bf x}_i:\ i=1, 2, \ldots \}$$
is a translative covering of $\mathbb{E}^n$ with density
$$\theta (\Im_1)\ge\theta^t(\sigma (K))=\theta^t(K).\eqno (17)$$
Furthermore, by (16) and (17),
$$\Im_2=\{rK_0+{\bf x}+{\bf x}_i:\ i=1, 2, \ldots \}$$
covers $\mathbb{E}^n$ with density $\theta (\Im_2)$, which satisfies
$$\theta (\Im_2)=r^n\theta^t(K_0)\ge \theta (\Im_1)\ge \theta^t(K).$$
The $f(K)=\theta^t(K)$ case of the theorem is proved. The other covering case
can be dealt with by similar arguments.

Assume that
$$\wp_0=\{rK+{\bf x}+{\bf x}_i:\ i=1, 2, \ldots \}$$
is a translative packing in $\mathbb{E}^n$ with the density
$$\delta (\wp_0)=\delta^t(K).$$
Then, by (16) it is easy to see that
$$\wp_1=\{ \sigma (K_0)+{\bf x}_i:\ i=1, 2, \ldots \}$$
is a translative packing in $\mathbb{E}^n$ with density
$$\delta (\wp_1)\le\delta^t(\sigma (K_0))=\delta^t(K_0).\eqno (18)$$
Furthermore, by (16) and (18),
$$\wp_2=\{K+{\bf x}_i:\ i=1, 2, \ldots \}$$
is a translative packing in $\mathbb{E}^n$ with density $\delta (\wp_2)$, which satisfies
$$\delta (\wp_2)=r^{-n}\delta^t(K)\le \delta (\wp_1)\le \delta^t(K_0).$$
Thus, we get
$$\delta^t(K)\le r^n\delta^t(K_0),$$
which proves the $f(K)=\delta^t(K)$ case of the theorem. The other packing case can be dealt with by similar arguments. \hfill{$\square$}

\medskip\noindent
{\bf Corollary 2.} {\it Let $f(K)$ to be $\theta^t(K)$, $\theta^l(K)$, $\delta^t(K)$ or $\delta^l(K)$. There are two positive constants $c(n)$ and $d(n)$, which depend on $n$, such that
$$|f(K_1)-f(K_2)|\le c(n)\cdot \| K_1, K_2\|$$
holds for any pair $K_1$, $K_2\in \mathcal{K}^n$, provided $\| K_1, K_2\|\le d(n)$. In particular, for all
$K\in \mathcal{K}^n$, we have}
$$f'(K)\le n.$$

\medskip\noindent
{\bf Remark 4.} When $n$ is small, both $c(n)$ and $d(n)$ can be taken precisely. For example, in the covering cases, one can take $c(2)=6$ and $d(2)={1\over 4}$, and take $c(3)=61$ and $d(3)={1\over 6}$; in the packing cases, one can take $c(2)=4$ and $d(2)={1\over 4}$, and take $c(3)=6$ and $d(3)={1\over 6}$.

\bigskip
Let $\phi^t(C)$ and $\phi^l(C)$ be the packing-covering constants defined in the previous section. For them we have the following results:

\medskip
\noindent
{\bf Theorem 4.} {\it Let $f(C)$ to be $\phi^t(C)$ or $\phi^l(C)$, and let $C_0$ be an $n$-dimensional centrally symmetric convex body. Then, for any $C\in \mathcal{C}^n$, we have}
$$f(C)\le e^{\| C, C_0\| }f(C_0).$$

\medskip\noindent
{\bf Proof.} As an example, we proceed to show the $f(C)=\phi^t(C)$ case. For convenience, we write
$$r=e^{\| C, C_0\|}.$$
By the definition of $\|\cdot\|$ we have
$$C\subseteq \sigma (C_0)\subseteq rC,\eqno(19)$$
where $\sigma $ is a suitable non-singular linear transformation from $\mathbb{E}^n$ to $\mathbb{E}^n$.
Furthermore, let $\phi (C,X)$ denote the ratio $\rho'/\rho $, where $\rho'$ is the smallest number such that $\rho'C+X$ is a covering of $\mathbb{E}^n$ and $\rho $ is the largest number such that $\rho C+X$ is a packing
in $\mathbb{E}^n$.

Let $X$ to be a discrete set in $\mathbb{E}^n$ such that
$$\phi (\sigma (C_0), X)={{\rho'}\over {\rho }}=\phi^t(C_0).\eqno(20)$$
In other words, $\rho'\sigma(C_0)+X$ is a covering of $\mathbb{E}^n$ and $\rho \sigma(C_0)+X$ is a packing in $\mathbb{E}^n$. Then, by (19), one can deduce that
$\rho'r C+X$ is a covering of $\mathbb{E}^n$, $\rho C+X$ is a packing in $\mathbb{E}^n$, and therefore, by (20),
$$\phi^t(C)\le \phi (C, X)\le {{\rho'r}\over \rho }=r\phi^t(C_0).$$
The translative case is proved. Clearly the lattice case can be shown by similar arguments. \hfill{$\square$}

\medskip\noindent
{\bf Corollary 3.} {\it Let $f(C)$ to be $\phi^t(C)$ or $\phi^l(C)$. There are two positive constants $c'(n)$ and $d'(n)$, which depend on $n$, such that
$$|f(C_1)-f(C_2)|\le c'(n)\cdot \| C_1, C_2\|$$
holds for any pair $C_1$, $C_2\in \mathcal{C}^n$, provided $\| C_1, C_2\|\le d'(n)$. In particular, for all
$C\in \mathcal{C}^n$, we have}
$$f'(C)\le 1.$$

\medskip\noindent
{\bf Remark 5.} When $n$ is small, both $c'(n)$ and $d'(n)$ can be taken precisely. It was proved by Zong \cite{zong08} and \cite{zong03} that $\phi^l(C)\le 2(2-\sqrt{2})$ holds for all two-dimensional centrally symmetric convex domains and $\phi^l(C)\le {7\over 4}$ for all three-dimensional centrally symmetric convex bodies.
Thus, one can take $c'(2)={5\over 2}$ and $d'(2)={1\over 2}$, and take $c'(3)={7\over 2}$ and $d'(3)={1\over 2}$.

\medskip
\noindent {\bf Problem 9.} Let $f(K)$ to be a good functional defined on $\{ \mathcal{K}^n, \|\cdot\|\}$, such as
$\delta^t(K)$, $\delta^l(K)$, $\delta^t(K)-\delta^l(K)$, $\theta^t(K)$, $\theta^l(K)$, $\theta^l(K)-\theta^t(K)$
and etc. To determine the value of
$$s(f'):=\max_{K\in \mathcal{K}^n}f'(K).$$

\medskip
\noindent {\bf Conjecture 1.} Let $f(K)$ to be a good functional defined on $\{ \mathcal{K}^n, \|\cdot\|\}$ such that
$f'(K)$ exists at every $K\in \mathcal{K}^n$. Then, for any pair of $K_1$, $K_2\in \mathcal{K}^n$, we have
$$|f(K_1)-f(K_2)|\le s(f')\cdot \| K_1, K_2\|.$$

\bigskip
\section*{5. Nets and Integral Sums}

\medskip
Let $\beta $ be a positive number and let $\mathcal{X}$ be a $\beta $-net in $\{ \mathcal{K}^n, \| \cdot \|\}$ with $\ell (n, \beta )$ elements. If $f(K)$ is a continuous functional defined on $\{ \mathcal{K}^n, \| \cdot \|\}$ such that
$$| f(K_1)-f(K_2)|\le c\cdot \| K_1, K_2\|$$
holds for any pair of convex bodies $K_1$, $K_2\in \mathcal{K}^n$ with some suitable constant $c$ provided $\| K_1, K_2\|\le \beta, $ then
$$\min_{K_i\in \mathcal{X}}f(K_i)-c\beta \le f(K)\le \max_{K_i\in \mathcal{X}}f(K_i)+ c\beta $$
holds for all $K\in \mathcal{K}^n$. Thus, by checking the values of $f(K_i)$ at $\ell (n, \beta )$ convex bodies, one
can estimate both
$$\max_{K\in \mathcal{K}^n} f(K)$$
and
$$\min_{K\in \mathcal{K}^n} f(K).$$
Clearly, all $\delta^t(K)$, $\delta^l(K)$, $\delta^t(K)-\delta^l(K)$, $\theta^t(K)$, $\theta^l(K)$ and $\theta^l(K)-\theta^t(K)$ can be dealt with in this way. Of course, $\ell (n, \beta )$ can be very huge, good $\beta $-net is hard to be constructed, and the values of $f(K_i)$ for particular $K_i$ are difficult to be determined. Nevertheless, this strategy does provide a theoretic mean to deal with many basic problems such as Problems 1, 2, 4, 5 and 8.

\medskip
If $\{ \mathcal{K}^n, \|\cdot\|\}$ has a good geometric measure and $f(K)$ is a good functional defined on
$\{ \mathcal{K}^n, \|\cdot\|\}$, one can define an integral, by which one can understand the average behavior of
$f(K)$ in $\{ \mathcal{K}^n, \|\cdot\|\}$. Unfortunately, as mentioned in Section 2, such good measure does not exist. However, we can introduce an integral sum on a net.

Let $\mathcal{X}=\left\{ K_1, K_2, \ldots , K_{\ell (n,\beta )}\right\}$ be a $\beta$-net in $\{ \mathcal{K}^n, \|\cdot\|\}$ with $\ell (n, \beta )$ elements and let $f(K)$ be a continuous functional defined on $\{ \mathcal{K}^n, \|\cdot\|\}$. Then, we define the following integral sum
$$\varpi (\mathcal{X}, f)={1\over {\ell (n, \beta )}}\sum_{K_i\in \mathcal{X}}f(K_i).$$
Clearly, this sum depends on the particular net $\mathcal{X}$. However, it reflects some measure theoretic property of $f(K)$. Of course, it also make sense to define and study sums based on ball packing systems in $\{ \mathcal{K}^n, \|\cdot\|\}$. We end this article with the following problem:

\medskip
\noindent
{\bf Problem 10.} For every continuous functional $f(K)$ in $\{ \mathcal{K}^n, \|\cdot\|\}$, does
$$\lim_{\beta\to 0}\varpi (\mathcal {X}, f)$$
always exist?

\vspace{1cm}\noindent
{\bf Acknowledgements.} This work is supported by 973 Programs 2013CB834201 and 2011CB302401, the National Science Foundation of China (No.11071003), and the Chang Jiang Scholars Program of China. I am grateful to Professor Senlin Wu for helpful comments.

\bibliographystyle{amsplain}

\end{document}